\documentclass[11pt]{amsart}
\usepackage{amsfonts, amssymb,enumerate}
\usepackage{url}

\newtheorem{thm}{Theorem}[section]
\newtheorem{prop}[thm]{Proposition}
\newtheorem{lem}[thm]{Lemma}
\newtheorem{cor}[thm]{Corollary}

\theoremstyle{definition}

\newtheorem{remark}[thm]{Remark}
\newtheorem{example}[thm]{Example}

\begin{document}

\newcommand\Aut{\operatorname{Aut}}
\newcommand\Ass{\operatorname{Ass}}
\newcommand\trdeg{\operatorname{trdeg}}
\newcommand\diag{\operatorname{diag}}
\newcommand\ed{\operatorname{ed}}
\newcommand\cd{\operatorname{cd}}
\newcommand\im{\operatorname{im}}
\newcommand\id{\operatorname{id}}
\newcommand\wt{\operatorname{wt}}
\newcommand\ind{\operatorname{ind}}
\newcommand\rk{\operatorname{rk}}
\newcommand\sign{\operatorname{sign}}
\newcommand\Spec{\operatorname{Spec}}
\newcommand\Spin{\operatorname{Spin}}
\newcommand\SO{\operatorname{SO}}
\newcommand\Mat{\operatorname{M}}
\newcommand\Ima{\operatorname{Im}}
\newcommand\Char{\operatorname{char}}
\newcommand\Gal{\operatorname{Gal}}
\newcommand\rank{\operatorname{rank}\,}
\newcommand\Ker{\operatorname{Ker}}

\newcommand\GL{\operatorname{GL}}
\newcommand\SL{\operatorname{SL}}
\newcommand\PGL{\operatorname{PGL}}
\newcommand\PGLn{\operatorname{PGL_n}}
\newcommand\Sym{\operatorname{S}}
\newcommand\Alt{\operatorname{A}}
\newcommand\Symn{\operatorname{S}_n}
\newcommand\Altn{\operatorname{A}_n}
\newcommand\Span{\operatorname{Span}}
\newcommand\bbC{\mathbb{C}}
\newcommand\ZZ{\mathbb{Z}}
\newcommand\bbZ{\mathbb{Z}}
\newcommand\bbG{\mathbb{G}}
\newcommand\bbQ{\mathbb{Q}}
\newcommand\bbF{\mathbb{F}}

\title[Essential dimension]{Some consequences
of the Karpenko-Merkurjev theorem}
\author[Aurel Meyer]{Aurel Meyer$^\dagger$}
\thanks{$^\dagger$ Aurel Meyer is partially supported by
a University Graduate Fellowship at the University of British Columbia}
\author[Zinovy Reichstein]{Zinovy Reichstein$^{\dagger \dagger}$}
\thanks{$^{\dagger \dagger}$ Z. Reichstein is partially supported by
NSERC Discovery and Accelerator Supplement grants}

\address{Department of Mathematics, University of British Columbia,
Vancouver, BC V6T 1Z2, Canada}


\begin{abstract} 
We use a recent theorem of N. A. Karpenko and 
A. S. Merkurjev to settle several  
questions in the theory of essential 
dimension. 
\end{abstract}

\maketitle

\section{Introduction}

N. Karpenko and A. Merkurjev~\cite{km} recently proved the following 
formula for the essential dimension of a finite $p$-group.

 \begin{thm} \label{thm.km} Let $G$ be a finite $p$-group and
 $k$ be a field of characteristic $\ne p$ containing 
 a primitive $p$th root of unity.
 Then $\ed(G; p) = \ed(G)$ = the minimal value of 
 $\dim(V)$, where the minimum is taked over all faithful linear 
 $k$-representations $G \to \GL(V)$. 
 \end{thm}

The purpose of this paper is to explore some of the consequences 
of this theorem.  We refer the reader 
to~\cite{br} or \cite[Chapter 8]{jly}
for background material on the essential dimension of a finite group,
\cite{bf} or~\cite{brv2} for the notion of essential dimension 
of a functor, and \cite{merkurjev} for a detailed discussion 
of essential dimension at a prime $p$.
As usual, if the reference to $k$ is clear from the context, we
will sometimes write $\ed$ in place of $\ed_k$.

The following notation will be used throughout.

For a finite group $H$, we will denote the intersection of 
the kernels of all multiplicative characters $\chi \colon H \to k^*$
by $H'$. 
In particular, if $k$ contains an $e$th root of unity,
where $e$ is the exponent of $H$, then $H' = [H, H]$. 

Given a $p$-group $G$, set $C(G)$ to be the center of $G$, $C(G)_p$ 
to be the $p$-torsion subgroup of $C(G)$. We will view $C(G)_p$ and its
subgroups as $\mathbb F_p$-vector spaces, and write ``$\dim_{\mathbb F_p}$" (or 
simply ``$\dim$") for their dimensions.  We further set
\begin{equation}\label{e.Ki}
K_i := \bigcap_{[G:H] = p^i}  H' \quad 
\text{and} \quad C_i := K_i \cap C(G)_p \, . 
\end{equation}
for every $i \ge 0$,
$K_{-1} := G$ and $C_{-1}:= K_{-1} \cap C(G)_p = C(G)_p$.

Our first main result is following theorem. Part (b) may be viewed 
as a variant of Theorem~\ref{thm.km}.

\begin{thm} \label{thm1} Let $G$ be a $p$-group, 
$k$ be a field of characteristic $\ne p$ containing 
a primitive $p$th root of unity, and $\rho \colon G \hookrightarrow \GL(V)$ 
be a faithful linear representation of $G$.  Then

\smallskip
(a) $\rho$ has minimal dimension among the faithful linear representations 
of $G$ defined over $k$ if and only if for every $i \ge 0$ 
the irreducible decomposition 
of $\rho$ has exactly \[ \dim_{\bbF_p} (C_{i-1}) - \dim_{\bbF_p} (C_i) \]
irreducible components of dimension $p^i$, each with multiplicity $1$.

\smallskip
(b) $\ed(G; p) = \ed(G) = \sum_{i = 0}^{\infty} 
(\dim_{\bbF_p} \, C_{i-1} - \dim_{\bbF_p} \, C_i) p^i$. 
\end{thm}

Note that $K_i = C_i = \{ 1 \}$ for large $i$ (say, if $p^i \ge |G|$),
so only finitely many terms in the above infinite sum are non-zero. 

We will prove Theorem~\ref{thm1} in section~\ref{sect.thm1};
the rest of the paper will be devoted to its applications.
The main results we will obtain are summarised below.

\subsection*{Classification of $p$-groups of essential dimension $\le p$.}

\begin{thm} \label{thm.ed=p}
Let $p$ be a prime, $k$ be a field of characteristic $\ne p$ containing 
a primitive $p$th root of unity and $G$ be a finite $p$-group 
such that $G' \ne \{ 1 \}$. Then the following conditions are equivalent.

\smallskip
(a) $\ed_k(G) \le p$,

\smallskip
(b) $\ed_k(G) = p$,

\smallskip
(c) The center $C(G)$ is cyclic and $G$ has a subgroup 
$A$ of index $p$ such that $A' = \{ 1 \}$.
\end{thm}

Note that the assumption that $G' \ne \{ 1 \}$ is harmless.
Indeed, if $G' = \{ 1 \}$ then by Theorem~\ref{thm1}(b) 
$\ed(G) = \rank(G)$; cf. also~\cite[Theorem 6.1]{br} 
or~\cite[section 3]{bf}. 

\subsection*{Essential dimension of $p$-groups of nilpotency class $2$.}

\begin{thm}\label{thm.p-groups}   
Let $G$ be a $p$-group of exponent $e$ and
$k$ be a field of characteristic $\ne p$ containing
a primitive $e$-th root of unity. Suppose the commutator
subgroup $[G, G]$ is central in $G$. Then

\smallskip
(a) $\ed_k(G; p) = \ed_k(G) \le 
\rank C(G) + \rank [G, G](p^{\lfloor m/2\rfloor}-1)$, 
where $p^m$ is the order of $G/C(G)$.

\smallskip
(b) Moreover, if $[G, G]$ is cyclic then $|G/C(G)|$ is a complete
square and equality holds in (a). That is, in this case
\[ \ed_k(G; p) = \ed_k(G) =\sqrt{|G/C(G)|} + \rank \, C(G) - 1 \, . \]
\end{thm}                                

\subsection*{Essential dimension of a quotient group}

C. U. Jensen, A. Ledet and N. Yui asked if $\ed G \ge \ed(G/N)$ for every
finite group $G$ and normal subgroup $N \triangleleft G$; 
see~\cite[p. 204]{jly}. The following theorem shows that
this inequality is false in general.

\begin{thm} \label{thm.jly}  
Let $p$ be a prime and $k$ be a field containing a primitive $p$th 
rooth of unity.
For every real number $\lambda > 0$ there exists 
a finite $p$-group $G$ and a central subgroup $H$ of $G$ 
such that $\ed_k(G/H) > \lambda \ed G$.
\end{thm} 

\subsection*{Essential dimension of $\SL_n(\bbZ)$}

G. Favi and M. Florence~\cite{ff} showed that
$\ed(\GL_n(\bbZ)) = n$ for every $n \ge 1$ and $\ed(\SL_n(\bbZ)) = n-1$ 
for every odd $n$. For details, including the definitions 
of $\ed(\GL_n(\bbZ))$ and
$\ed(\SL_n(\bbZ))$, see Section~\ref{sect.ff}. For even $n$ Favi and 
Florence showed that $\ed(\SL_n(\bbZ)) = n-1$ or $n$ and left the exact
value of $\SL_n(\bbZ)$ as an open question. In this paper
we will answer this question as follows. 

\begin{thm} \label{thm.ff1}
Suppose $k$ is a field of characteristic $\ne 2$. Then
\[ \ed_k (\SL_n(\bbZ);2) = 
\ed_k (\SL_n(\bbZ)) = \begin{cases} \text{$n -1$, if $n$ is odd,} \\
\text{$n$, if $n$ is even} \end{cases}
\]
for any $n \ge 3$.
\end{thm}

\subsection*{Acknowledgement}
Theorems~\ref{thm.p-groups}(b) and~\ref{thm.jly}(b) first 
appeared in the unpublished preprint~\cite{brv1} by P. Brosnan, the second
author and A. Vistoli. We thank P. Brosnan and A. Vistoli for 
allowing us to include them in this paper.  Theorem~\ref{thm.p-groups}(b) 
was, in fact, a precursor to Theorem~\ref{thm.km}; the techniques 
used in~\cite{brv1} were subsequently strengthened and refined 
by Karpenko and Merkurjev~\cite{km} to prove Theorem~\ref{thm.km}. 
The proof of Theorem~\ref{thm.p-groups}(b) in Section~\ref{sect.brv}
may thus be viewed as a result of reverse engineering. We include 
it here because it naturally fits into the framework 
of this paper, because Theorem~\ref{thm.p-groups}(b) is used 
in a crucial way in~\cite{brv2}, and because a proof of this result 
has not previously appeared in print.

\section{Proof of Theorem~\ref{thm1}}
\label{sect.thm1}

Throughout this section $G$ will denote a $p$-group, 
and $C = C(G)_p$ will denote the $p$-torsion subgroup 
of the center of $G$. We will use the notations introduced 
in and just before the statement of Theorem~\ref{thm1}.
In particular,
\[ K_{-1} = G \supset K_0 \supset K_1 \supset K_2 \supset \dots \]
will be the descending sequence of normal subgroups of $G$ defined in
\eqref{e.Ki} and $C_i = C \cap K_i$. 
We will repeatedly use the well-known fact that
\begin{equation} \label{fact}
\text{A normal subgroup $N$ of $G$ is trivial if and only if $N \cap C$ is
trivial.}
\end{equation}

We begin with three elementary lemmas.

\begin{lem}\label{l.Ki} 
$K_i=\bigcap_{\dim(\rho) \le p^i} \, \ker(\rho)$, where the intersection is
taken over all irreducible representations $\rho$ of $G$ of dimension $\le p^i$.
\end{lem}

\begin{proof}
Let $j\le i$.
Recall that every irreducible representation $\rho$ of $G$ of 
dimension $p^j$ is induced from a $1$-dimensional representation 
$\chi$ of a subgroup $H \subset G$ of index $p^j$; cf.~\cite[8.16]{serre-rep}.
Thus $\ker(\rho)=\ker(\ind_H^G\chi)=\bigcap_{g\in G}g\ker(\chi)g^{-1}$, 
and since each $g\ker(\chi)g^{-1}$ contains $(gHg^{-1})^\prime$, 
we see that $\ker(\rho)\supset K_j\supset K_i$.
The opposite inclusion is proved in a similar manner.
\end{proof}

\begin{lem}\label{l.rank}
Let $G$ be a finite group over a field $k$ that contains $p$th roots
of unity. Let $C$ be a central subgroup of exponent $p$
and $\rho \colon G \to \GL(V)$ be an irreducible representation
of $G$. Then 

\smallskip
(a) $\rho(C)$ consists of scalar matrices.
In other words, the restriction of $\rho$ to $C$ decomposes
as $\chi \oplus \ldots \oplus \chi$ ($\dim(V)$ times),
for some multiplicative character $\chi \colon C \to \bbG_m$.
We will refer to $\chi$ as the {\em character associated to} $\rho$.

\smallskip
(b) $C_i=\bigcap_{\dim(\rho) \le p^i} \,  \ker(\chi_{\rho})$, where 
the intersection is taken over all irreducible 
$G$-representations $\rho$ of dimension $\le p^i$ and
$\chi_{\rho} \colon C \to \bbG_m$ denotes the character associated 
to $\rho$. 

In particular, if $\dim(\rho) \le p^i$ then the associated 
character $\chi$ of $\rho$ vanishes on $C_i$.
\end{lem}

\begin{proof}
(a) Let $V = V_1 \oplus \dots \oplus V_m$ be an irreducible decomposition
of $V$ as a direct sum of character spaces for $C$. That is, $C$ acts
on $V_i$ by a multiplicative character $\chi_i$, where $\chi_1, \dots,
\chi_m$ are distinct. Since $C$ is central, each $V_i$ is $G$-invariant.
Since we are assuming that the representation of $G$ on $V$ 
is irreducible, this implies that $m = 1$, as claimed.

\smallskip
(b) By Lemma~\ref{l.Ki} 
\[ C_i = C \cap \bigcap_{\dim(\rho) \le p^i} \, \ker(\rho) = 
\bigcap_{\dim(\rho) \le p^i} \, (C \cap \ker(\rho)) = 
\bigcap_{\dim(\rho) \le p^i} \, \ker(\chi_{\rho}) \, . \]
\end{proof}

\begin{lem} \label{lem2.3} Let $G$ be a $p$-group and
$\rho = \rho_1 \oplus \ldots \oplus \rho_m$ be the direct sum of
the irreducible representations $\rho_i \colon G \to \GL(V_i)$.
Let $\chi_i := \chi_{\rho_i} \colon C \to \bbG_m$ be the character
associated to $\rho_i$.

\smallskip
(a) $\rho$ is faithful if and only if $\chi_1, \dots, \chi_m$ span $C^*$
as an $\mathbb F_p$-vector space.

\smallskip
(b) Moreover, if $\rho$ is of minimal dimension among the faithful
representations of $G$ then  $\chi_1, \dots, \chi_m$ form 
an $\mathbb F_p$-basis of $C^*$.
\end{lem}

\begin{proof} (a) By~\eqref{fact},
$\Ker(\rho)$ is trivial if and only if 
$\Ker(\rho) \cap C = \cap_{i = 1}^m \Ker(\chi_i)$ is trivial.
On the other hand, 
$\cap_{i = 1}^m \Ker(\chi_i)$ is trivial if and only if
$\chi_1, \dots, \chi_m$ span $C^*$.

\smallskip
(b) Assume the contrary, say $\chi_m$ is a linear combination of
$\chi_1, \dots, \chi_{m-1}$. Then part (a) tells us that
$\rho_1 \oplus \ldots \oplus\rho_{m-1}$ is a faithful representation of $G$, 
contradicting the minimality of $\dim(\rho)$.
\end{proof}

We are now ready to proceed with the proof of Theorem~\ref{thm1}.
Part (b) is an immediate consequence of part (a) and Theorem~\ref{thm.km}.
We will thus focus on proving part (a). 
In the sequel for each $i \ge 0$ we will set 
\[ \delta_i := (\dim \, C_{i-1} - \dim \, C_i) \]
and
\[ \Delta_i : = \delta_0 + \delta_1 + \dots + \delta_i = 
\dim(C) - \dim(C_i) \, , \]
where the last equality follows from $C_{-1} = C$.

Our proof will proceed in two steps. In Step 1 we will construct
a faithful representation $\mu$ of $G$ such that for every $i \ge 0$
exactly $\delta_i$ irreducible components of $\mu$ have dimension $p^i$.
In Step 2 we will show that $\dim(\rho) \ge \dim(\mu)$ 
for any other faithful representation $\rho$ of $G$, and moreover 
equality holds if and only if for every $i \ge 0$ 
$\rho$ has exactly $\delta_i$ irreducible components 
of dimension $p^i$.

\smallskip
{\bf Step 1:}
We begin by constructing $\mu$.  By definition, 
\[ C = C_{-1} \supset C_{0} \supset C_1 \supset \dots  \, ,  \] 
where the inclusions are not necessarily strict.
Dualizing this flag of $\mathbb F_p$-vector spaces, we 
obtain a flag 
\[ (0) = (C^*)_{-1} \subset (C^*)_0 \subset (C^*)_1 \subset \dots
\] 
of $\mathbb F_p$-subspaces of $C^*$, where
\[ \text{$(C^*)_i := \{ \chi \in C^* \, | \, \chi$ is trivial 
on $C_i \} \simeq (C/C_i)^*$.} \]
Let $\Ass(C) \subset C^*$ be the set of characters of $C$ associated 
to irreducible representations of $G$, and let $\Ass_i(C)$ be the set
of characters associated to irreducible representations of 
dimension $p^i$. Lemma~\ref{l.rank}(b) tells us that
\[ \text{$\Ass_0(C) \cup \Ass_1(C) \cup \dots \cup \Ass_i(C)$ spans $(C^*)_i$} \]
for every $i \ge 0$.  Hence, we can choose a basis 
$\chi_1, \dots, \chi_{\Delta_0}$ of $(C^*)_0$ from $\Ass_0(C)$, then
complete it to a basis 
$\chi_1, \dots, \chi_{\Delta_1}$ of $(C^*)_1$ by choosing
the last $\Delta_1 - \Delta_0$ characters from $\Ass_1(C)$, 
then complete this basis of $(C^*)_1$ to a basis of $(C^*)_2$ by choosing
$\Delta_2 - \Delta_1$ additional characters from $\Ass_2(C)$, 
etc. 
We stop when $C_i = (0)$, i.e., $\Delta_i = \dim(C)$.  

By the definition of $\Ass_i(C)$, each $\chi_j$ is the associated character 
of some irreducible representation $\mu_j$ of $G$. By our construction
\[ \mu = \mu_1 \oplus \dots \oplus \mu_{\dim(C)} \, , \]
has the desired properties. 
Indeed, since 
$\chi_1, \dots, \chi_{\dim(C)}$ form a basis of $C^*$, Lemma~\ref{lem2.3}
tells us that $\mu$ is faithful. 
On the other hand, by 
our construction exactly
\[ \delta_i - \delta_{i-1} = \dim(C_i^*) - \dim(C_{i-1}^*) =
\dim(C_{i-1}) - \dim(C_i) \] 
of the characters $\chi_1, \dots, \chi_c$ come from $\Ass_i(C)$. 
Equivalently, exactly
$\dim(C_{i-1}) - \dim(C)$ of the irreducible representations
$\mu_1, \dots, \mu_c$ are of dimension $p^i$. 

\medskip
{\bf Step 2:} Let $\rho \colon G \to \GL(V)$ be a faithful linear 
representation of $G$ of the smallest possible dimension, 
\[ \rho = \rho_1 \oplus \ldots \oplus \rho_c \] 
be its irreducible decomposition, and $\chi_i \colon C \to \bbG_m$ 
be the character associated to $\rho_i$. By Lemma~\ref{lem2.3}(b), 
$\chi_1, \dots, \chi_c$ form a basis of $C^*$. 
In particular, $c = \dim(C)$ 
and at most $\dim(C) - \dim(C_i)$ of the characters 
$\chi_1, \dots, \chi_c$
can vanish on $C_i$. On the other hand, by
Lemma~\ref{l.rank}(b) every representation of dimension
$\le p^i$ vanishes on $C_i$.  Thus if 
exactly $d_i$ of the irreducible representations 
$\rho_1, \dots, \rho_c$ have dimension $p^i$ then
\[ d_0 + d_1 + d_2 + \ldots + d_i \le \dim(C) - \dim(C_i) 
\]
for every $i \ge 0$. For $i \ge 0$, set
$D_i := d_0 + \dots + d_i$ = number of representations of dimension $\le p^i$
among $\rho_1, \dots, \rho_c$. 
We can now write the above inequality as
\begin{equation} \label{e.inequality}
\text{$D_i \le \Delta_i$ for every $i \ge 0$.} 
\end{equation}
Our goal is to show that
$\dim(\rho) \ge \dim(\mu)$ and that equality holds if and only if
exactly $\delta_i$ of the irreducible representations 
$\rho_1, \dots, \rho_{\dim(C)}$ have dimension $p^i$.
The last condition translates into $d_i = \delta_i$ 
for every $i \ge 0$, which is, in turn equivalent to
$D_i = \Delta_i$ for every $i \ge 0$. 

Indeed, setting $D_{- 1} := 0$ and $\Delta_{-1} := 0$, we have,
\begin{eqnarray*}
\dim(\rho)  - \dim(\mu) = \sum_{i = 0}^{\infty} (d_i - \delta_i) p^i =  
\sum_{i = 0}^{\infty} (D_i - \Delta_i)p^i - 
\sum_{i = 0}^{\infty} (D_{i-1} - \Delta_{i-1})p^i \\ 
= \sum_{i = 0}^{\infty} (D_i - \Delta_i)(p^i - p^{i + 1}) \ge 0 \, ,    
\end{eqnarray*}
where the last inequality follows from~\eqref{e.inequality}. Moreover,
equality holds if and only if $D_i = \Delta_i$ for every $i \ge 0$, as claimed.
This completes the proof of Step 2 and thus of Theorem~\ref{thm1}.
\qed

\section{Proof of Theorem~\ref{thm.ed=p}}

Since $K_0 = G'$ is a non-trivial normal subgroup of $G$, we see
that $K_0 \cap C(G)$ and thus $C_0 = K_0 \cap C(G)_p$ is non-trivial. 
This means that in the summation formula of Theorem~\ref{thm1}(b) at
least one of the terms 
\[ (\dim_{\bbF_p}(C_{i-1}) - \dim_{\bbF_p}(C_i)) p^i \]
with $i \ge 1$ will be non-zero. Hence,
$\ed(G) \ge p$; this shows that (a) and (b) are equivalent.
Moreover, equality holds 
if and only if (i) $\dim_{\bbF_p}(C_{-1}) = 1$, (ii) $\dim_{\bbF_p}(C_0) = 1$ 
and (iii) $C_1$ is trivial. 
It remains to show that (i), (ii) and (iii) are equivalent to (c).

Since $C_{-1} = C(G)_p$, (i) is equivalent to $C(G)$ being cyclic. 

Now recall that we are assuming $K_0 = G' \ne \{ 1 \}$. By~\eqref{fact}
this is equivalent to
$C_0 = K_0 \cap C(G)_p \ne \{ 1 \}$. Since $C_0 \subset C_{-1}$ 
has dimension at most 1, we see that (ii) follows from (i).

Finally, (iii) means that
\begin{equation} \label{e.K_1}
K_1 = \bigcap_{[G:H] = p} \, H'
\end{equation}
intersects $C(G)_p$ trivially.  Since $K_1$ is a normal subgroup 
of $G$, \eqref{fact} tells us that (iii) holds if and only 
if $K_1 = \{ 1 \}$. 

It thus remains to show that $K_1 = \{ 1 \}$ if and only
if $H' = \{ 1 \}$ for some subgroup $H$ of $G$ of index $p$.
One direction is obvious: if $H' = \{ 1 \}$ for some $H$ of index $p$
then the intersection~\eqref{e.K_1} is trivial.
To prove the converse, assume the contrary: the intersection~\eqref{e.K_1} 
is trivial but $H' \ne \{ 1 \}$ for every subgroup $H$ of index $p$.
Since every such $H$ is normal in $G$ (and so is $H'$),~\eqref{fact}
tells us that that $H' \ne \{ 1 \}$ if and only if 
$H' \cap C(G) \ne \{ 1 \}$. Since $C(G)$ is cyclic, the latter 
condition is equivalent to $C(G)_p \subset H'$. Thus
\[ C(G)_p \subset K_1 = \bigcap_{[G:H] = p} \,  H' \, , \]
contradicting our assumption that $K_1 \ne \{ 1 \}$.
\qed

\section{Proof of Theorems~\ref{thm.p-groups} and ~\ref{thm.jly}}
\label{sect.brv}
\begin{proof}[Proof of Theorem~\ref{thm.p-groups}]
Since the commutator $K_0=[G,G]$ is central, $C_0=K_0\cap C(G)_p$ 
is of dimension $\rank [G,G]$ and the $p^0$ term in 
the formula of Theorem~\ref{thm1} is $(\rank C(G)-\rank [G,G])$.

Let $Q=G/C(G)$ which is abelian by assumption. 
Let $h_1,...,h_s$ be generators of $[G,G]$ where 
$s=\rank [G,G]$, so that 
\[ [G,G]=\bbZ/p^{e_1}h_1\oplus
\dots \oplus \bbZ/p^{e_1}h_1 \, , \]
written additively.  For $g_1,g_2\in G$ the commutator can then 
be expressed as \[ [g_1,g_2]=\beta_1(g_1, g_2) h_1+ 
\ldots +\beta_s(g_1, g_2) h_s \, . \]
Note that each $\beta_i(g_1, g_2)$ depends on $g_1,g_2$ only
modulo the center $C(G)$. Thus each $\beta_i$ descends to a
skew-symmetric bilinear form
\[
Q\times Q\rightarrow \bbZ/p^{e_i} 
\]
which, by a slight abuse of notation, we will continue to denote 
by $\beta_i$. 
Let $p^m$ be the order of $Q$. For each form $\beta_i$
there is an isotropic subgroup $Q_i$ of $Q$ of order 
at least $p^{\lfloor (m+1)/2 \rfloor}$ (or equivalently, of
index at most $p^{\lfloor m/2 \rfloor}$ in $Q$); see \cite[Corollary 3]{at}.
Pulling these isotropic subgroups back to $G$, we obtain 
subgroups $G_1, \dots, G_s$ of $G$ of index $\le p^{\lfloor m/2 \rfloor}$ 
with the property that $G_i' = [G_i, G_i]$ lies in the subgroup of 
$C(G)$ generated by $h_1, \dots, h_{i-1}, h_{i+1}, \dots, h_s$.
In particular, $G_i' \cap \dots \cap G_s' = \{ 1 \}$.
Thus, all $K_i$ (and hence, all $C_i$) in \eqref{e.Ki} are trivial 
for $i\ge \lfloor m/2 \rfloor$, and Theorem~\ref{thm1} tells us that
\begin{eqnarray*} \ed(G) =  \dim_{\bbF_p}  \, C_{-1} - \dim_{\bbF_p}\, C_0 +
\sum_{i = 1}^{\lfloor m/2 \rfloor} 
(\dim_{\bbF_p} \, C_{i-1} - \dim_{\bbF_p} \, C_i)  p^i
\le  \\
 \dim_{\bbF_p}  \, C_{-1} - \dim_{\bbF_p} \, C_0 +
\sum_{i = 1}^{\lfloor m/2 \rfloor} 
(\dim_{\bbF_p} \, C_{i-1} - \dim_{\bbF_p} \, C_i)  \cdot p^{\lfloor m/2 \rfloor} = \\ 
\rank \, C(G)  +  p^{\lfloor m/2 \rfloor} (\rank \, [G, G] - 1) \, . 
\end{eqnarray*}

(b) In general, the skew-symmetric bilinear forms $\beta_i$ 
may be degenerate. However, 
if $[G,G]$ is cyclic, i.e., $s = 1$, then we have only one form, 
$\beta_1$, which is easily seen to be non-degenerate.  
For notational simplicity, we will write $\beta$ instead of $\beta_1$. 
To see that $\beta$ is non-degenerate,
suppose $\overline{g} := g$ (modulo $C(G)$) 
lies in the kernel of $\beta$ for some $g \in G$. Then by definition 
\[ \beta(g, g_1) = g g_1 g^{-1} g_1^{-1} = 1 \]
for every $g_1 \in G$. Hence, $g$ is central in $G$,
i.e., $\overline{g} = 1$ in $Q = G/C(G)$, as claimed.
     
We conclude that the order of $Q=G/C(G)$ is a perfect square, say $p^{2i}$, 
and $Q$ contains a maximal isotropic subgroup $I \subset Q$ of order 
$p^i = \sqrt{|G/C(G)|}$; see \cite[Corollary 4]{at}. 
The preimage of $I$ in $G$ is a maximal abelian subgroup of 
index $p^i$.  Consequently, $K_0= [G,G], K_1 , \ldots, 
K_{i-1}$ are all of rank $1$ and $K_i$ is trivial, 
where $p^i=\sqrt{|G/C(G)|}$.
Moreover, since all of these groups lie 
in $[G, G]$ and hence, are central, we have $C_i = (K_i)_p$ and thus
\[ \text{$\dim_{\mathbb F_p}(C_0) = \dim_{\mathbb F_p}(C_1) = 
\ldots = \dim_{\mathbb F_p}(C_{i-1}) = 1$ and
$\dim_{\mathbb F_p}(C_i) = 0$} \, . \] 
Specializing the formula of Theorem~\ref{thm.p-groups} to this
situation, we obtain part (b).
\end{proof}

\begin{example} \label{ex.extraspecial}
Recall that a $p$-group $G$ is called \emph{extra-special} if its
center $C$ is cyclic of order $p$, and the quotient $G/C$ is elementary
abelian. The order of an extra special $p$-group $G$ is an odd power
of $p$; the exponent of $G$ is either $p$ or $p^{2}$;
cf.~\cite[III. 13]{huppert}.  Note that every non-abelian
group of order $p^{3}$ is extra-special. For extra-special
$p$-groups Theorem~\ref{thm.p-groups}(b) reduces to the following.

{\em Let $G$ be an extra-special $p$-group of order $p^{2m+1}$.
Assume that the characteristic of $k$ is different from $p$,
that $\zeta_{p} \in k$, and $\zeta_{p^{2}} \in k$ if the exponent
of $G$ is $p^{2}$. Then $\ed G = p^{m}$.}   
\end{example}

\begin{proof}[Proof of Theorem~\ref{thm.jly}]
Let $\Gamma$ be a non-abelian group of order $p^3$.
The center of $\Gamma$ has order $p$; denote it by $C$. 
Since $\Gamma$ is extra-special, $\ed(\Gamma) = p$. (This 
also follows from Theorem~\ref{thm.ed=p}.)

The center of $\Gamma^n = \Gamma \times \dots \times \Gamma$
($n$ times) is then isomorphic to $C^n$. Let $H_n$ be the subgroup 
of $C^n$ consisting of $n$-tuples $(c_1, \dots, c_n)$ such that
$c_1 \dots c_n = 1$.  Clearly
   \[
   \ed\Gamma^n \le n \cdot \ed(\Gamma) = np \, ;
   \]
see \cite[Lemma 4.1(b)]{br}. (In fact 
by \cite[Theorem 5.1]{km}, $\ed\Gamma^n = n \cdot \ed(\Gamma)$ 
but we shall not need this here.)

On the other hand, $\Gamma^n/H_{n}$, is easily seen to be extra-special
of order $p^{2n+1}$, so $\ed(\Gamma^n/H_{n}) = p^{n}$ by
Example~\ref{ex.extraspecial}. Setting $G = \Gamma^n$ and $H = H_n$, 
we see that the desired inequality $\ed(G/H) > \lambda \ed G$ 
holds for suitably large $n$.
\end{proof}

\section{Proof of Theorem~\ref{thm.ff1}}  
\label{sect.ff}

Recall that the essential dimension of the group $\GL_n(\bbZ)$ 
over a field $k$, or $\ed_k(\GL_n(\bbZ))$ for short, is defined 
as the essential dimension of this functor
\[ H^1( \ast , \GL_n(\bbZ))\colon
 K \to \text{\{$K$-isomorphism classes of $n$-dimensional $K$-tori\}} \, , \]
where $K/k$ is a field extension.  Similarly
$\ed_k(\SL_n(\bbZ))$ is defined as the essential dimension 
of the functor
\[ \begin{array}{r} H^1( \ast , \SL_n(\bbZ)) \colon
K \to \text{\{$K$-isomorphism classes of $n$-dimensional $K$-tori} \\
\text{with $\phi_T \subset \SL_n(\bbZ)$ \}} \, , \end{array} \]
where $\phi_T \colon \Gal(K) \to \GL_n(\bbZ)$ is the natural representation
of the Galois group of $K$ on the character lattice of $T$.
The essential dimensions $\ed_k(\GL_n(\bbZ); p)$ and
$\ed_k(\SL_n(\bbZ); p)$ are respectively the essential dimensions 
of the above functors at a prime $p$.

G. Favi and M. Florence~\cite{ff} showed that for $\Gamma = \GL_n(\bbZ)$ or
$\SL_n(\bbZ)$, 
\begin{equation}\label{e.edL}
\ed(\Gamma)=\max \{\ed(F)| F\mbox{ finite subgroup of }\Gamma\}.
\end{equation} 
 From this they deduced that
\[ \ed(\GL_n(\bbZ)) = n, \quad \text{and} \quad
\ed(\SL_n(\bbZ)) = \begin{cases} \text{$n -1$, if $n$ is odd,} \\
\text{$n -1$ or $n$, if $n$ is even.} \end{cases} \]
For details, see \cite[Theorem 5.4]{ff}.

Favi and Florence also proved that $\ed(\SL_2(\bbZ)) = 1$
if $k$ contains a primitive 12th root of unity and asked whether
$\ed(\SL_n(\bbZ)) = n-1$ or $n$, in the case where $n \ge 4$ is
even; see~\cite[Remark 5.5]{ff}.  In this section we will prove
Theorem~\ref{thm.ff1} which shows that the answer is always $n$.

A minor modification of the arguments in \cite{ff} shows 
that~\eqref{e.edL} holds also for essential dimension at a prime $p$:
\begin{equation}\label{e.edLp}
\ed(\Gamma;p)=\max \{\ed(F;p)| F\mbox{ a finite subgroup of }\Gamma\},
\end{equation} 
where $\Gamma = \GL_n(\bbZ)$ or $\SL_n(\bbZ)$.
The finite groups $F$ that Florence and Favi used to find 
the essential dimension of $\GL_n(\bbZ)$ and $\SL_n(\bbZ)$ ($n$ odd) 
are $(\bbZ/2 \bbZ)^n$ and  $(\bbZ/2 \bbZ)^{n-1}$ respectively.
Thus $\ed(\GL_n(\bbZ); 2) = \ed(\GL_n(\bbZ)) = n$ for every $n \ge 1$ 
and $\ed(\SL_n(\bbZ); 2) = \ed(\SL_n(\bbZ)) = n-1$ if $n$ is odd. 

Our proof of Theorem~\ref{thm.ff1} will rely on part (b) of
the following easy corollary of Theorem~\ref{thm1}.

\begin{cor} \label{cor.km3} 
Let $G$ be a finite $p$-group, and $k$ be a field of characteristic
$\ne p$, containing a primitive $p$th root of unity.

\smallskip
(a) If $C(G)_p \subset K_i$ then $\ed_k(G)$ is divisible by $p^{i+1}$.

\smallskip
(b) If $C(G)_p \subset G'$ then $\ed_k(G)$ is divisible by $p$.  

\smallskip
(c) If $C(G)_p \subset G^{(i)}$, where $G^{(i)}$ denotes the $i$th 
derived subgroup of $G$, then $\ed_k(G)$ is divisible by $p^{i}$.
\end{cor}

\begin{proof} 
(a) $C(G)_p \subset K_i$ implies $C_{-1} = C_{0} = \dots = C_{i}$.
Hence, in the formula of Theorem~\ref{thm1}(b)
the $p^0, p^1, \dots, p^{i}$ terms appear with 
coefficient $0$. All other terms are divisible by $p^{i+1}$,
and part (a) follows.

(b) is an immediate consequence of (a), since $K_0 = G'$.

(c) By ~\cite[Theorem V.18.6]{huppert} $G^{(i)}$ is contained
in the kernel of every $p^{i-1}$-dimensional representation
of $G$. Lemma~\ref{l.Ki} now tells us that $G^{(i)} \subset K_{i-1}$
and part (c) follows from part (a).
\end{proof}

\begin{proof}[Proof of Theorem~\ref{thm.ff1}]
We assume that $n = 2d \ge 4$ is even. 
To prove Theorem~\ref{thm.ff1} it suffices 
to find a finite 2-subgroup $F$ of $\SL_n(\bbZ)$ of essential 
dimension $n$.

Diagonal matrices and permutation matrices generate a subgroup
of $\GL_n(\bbZ)$ isomorphic to $\mu_2^n \rtimes \Symn$. 
The determinant function restricts to a homomorphism  
$\det \colon \mu_2^n \rtimes \Symn \to \mu_2$
sending $((\epsilon_1, \dots, \epsilon_n), \tau)) \in \mu_2^n \rtimes \Symn$
to the product $\epsilon_1 \epsilon_2 \cdot \ldots \cdot \epsilon_n \cdot \sign(\tau)$. Let $P_n$ be 
a Sylow $2$-subgroup of $\Symn$ and $F_n$ be the kernel of
\[ \det \colon \mu_2^n \rtimes P_n \to \mu_2 \, . \]
By construction $F_n$ is a finite $2$-group contained in $\SL_n(\bbZ)$.
Theorem~\ref{thm.ff1} is now a consequence of the following proposition.

\begin{prop} \label{prop.ff2} If $\Char(k) \ne 2$ then 
$\ed(F_{2d}) = 2d$ for any $d \ge 2$.
\end{prop}

To prove the proposition,
let \[ D_{2d} = \{ \diag(\epsilon_1, \dots, \epsilon_{2d}) \, | \, 
\text{each $\epsilon_i = \pm 1$ and
$\epsilon_1 \epsilon_2 \cdots \epsilon_{2d} = 1$} \} \]
be the subgroup of ``diagonal" matrices contained in $F_{2d}$.

Since $D_{2d} \simeq \mu_2^{2d-1}$ has essential dimension $2d-1$, 
we see that $\ed(F_{2d}) \ge \ed(D_{2d}) = 2d - 1$.  On the other 
hand the inclusion $F_{2d} \subset \SL_{2d}(\bbZ)$ gives rise
to a $2d$-dimensional representation of $F_{2d}$, which remains faithful
over any field $k$ of characteristic $\neq 2$. 
Hence, $\ed(F_{2d}) \le 2d$. We thus conclude that 
\begin{equation} \label{e.ff2}
\text{$\ed(F_{2d}) = 2d - 1$ or $2d$.}
\end{equation} 
Using elementary group theory, one easily 
checks that \begin{equation} \label{e.check}
C(F_{2d}) \subset [F_{2d}, F_{2d}] \subset F_{2d}' \, . 
\end{equation}
Thus $\ed(F_{2d})$ is even by Corollary~\ref{cor.km3}; 
\eqref{e.ff2} now tells us that $\ed(F_{2d}) = 2d$.
This completes the proof of Proposition~\ref{prop.ff2} and thus 
of Theorem~\ref{thm.ff1}.
\end{proof}

\begin{remark} \label{rem.d=1}
The assumption that $d \ge 2$ is essential in the proof 
of the inclusion~\eqref{e.check}. In fact, $F_2 \simeq \bbZ/4 \bbZ$,
so~\eqref{e.ff2} fails for $d = 1$.
\end{remark}

\begin{remark}  
Note that for any integers $m, n \ge 2$, $F_{m + n}$ contains 
the direct product $F_m \times F_n$.  Thus 
\[ \ed(F_{m + n}) \ge  \ed(F_m \times F_n) = \ed(F_m) + \ed(F_n) \, ,  \]
where the last equality follows from~\cite[Theorem 5.1]{km}. 
Thus Proposition~\ref{prop.ff2} only needs to be proved for $d = 2$ 
and $3$ (or equivalently, $n = 4$ and $6$); 
all other cases are easily deduced from these by applying
the above inequality recursively, with $m = 4$.  In particular,
the group-theoretic inclusion~\eqref{e.check} only needs to be 
checked for $d = 2$ and $3$. Somewhat to our surprise, this 
reduction does not appear to simplify the proof of 
Proposition~\ref{prop.ff2} presented above to any 
significant degree.
\end{remark}

\begin{remark}
It is interesting to note that while 
the value of $\ed_k(\SL_2(\bbZ))$ depends on the 
base field $k$ (see~\cite[Remark 5.5]{ff}),
for $n \ge 3$, the value of $\ed_k(\SL_n(\bbZ))$ does not
(as long as $\Char(k) \ne 2$).
\end{remark}

\end{document}